\documentclass[12pt]{article}%
\usepackage{amssymb}
\usepackage{amsfonts}
\usepackage{sw20bams}%
\usepackage{amsmath}%
\usepackage{graphicx}
\providecommand{\U}[1]{\protect\rule{.1in}{.1in}}
\begin{document}

\title{Maximum Queue Length for Traffic Light with Bernoulli Arrivals}
\author{Steven Finch}
\date{February 13, 2018}
\maketitle

\begin{abstract}
Cars arrive at an intersection with a stoplight, which is either red or green.
The cars all travel in the same direction, that is, we ignore cross-traffic \&
oncoming traffic. Assume that the intersection is initially empty. Assume
that, at every second, there is a probability $p$ that one new car will arrive
at the light, and the outcome is independent of past \& future. Let $\ell
\geq1$ be an integer. A red light lasts $\ell$ seconds; likewise for green. If
the light is red, no cars can leave the intersection. If the light is green,
cars will leave the intersection at a rate of one per second. Over a time
period of $n$ seconds, determine the (random) maximum queue length $M$ of cars
at the intersection. What is the distribution of $M$, as a function of
$(p,\ell,n)$? \ We answer this question for the special case $\ell=1$ and
introduce a conjecture for $\ell>1$.

\end{abstract}

\footnotetext{Copyright \copyright \ 2018 by Steven R. Finch. All rights
reserved.}Let $\ell\geq1$ be an integer and $p>0$ be a real number with
$p+q=1$, $p\leq q$. \ Let $X_{0}=0$ and $X_{1}$, $X_{2}$, \ldots, $X_{n}$ be a
sequence of independent random variables satisfying%
\[%
\begin{array}
[c]{ccccc}%
\mathbb{P}\left\{  X_{i}=1\right\}  =p, &  & \mathbb{P}\left\{  X_{i}%
=0\right\}  =q &  & \text{if }i\equiv1,2,\ldots,\ell\,\operatorname{mod}%
\,2\ell;
\end{array}
\]%
\[%
\begin{array}
[c]{ccccc}%
\mathbb{P}\left\{  X_{i}=0\right\}  =p, &  & \mathbb{P}\left\{  X_{i}%
=-1\right\}  =q &  & \text{if }i\equiv\ell+1,\ell+2,\ldots,2\ell
\,\operatorname{mod}\,2\ell
\end{array}
\]
for each $1\leq i\leq n$. \ Define $S_{0}=X_{0}$ and $S_{j}=\max\left\{
S_{j-1}+X_{j},0\right\}  $ for all $1\leq j\leq n$. \ The (non-simple)
reflected random walk $S_{0}$, $S_{1}$, $S_{2}$, \ldots, $S_{n}$ serves as a
primitive model for traffic queue congestion at a stoplight, where arrivals
follow a Bernoulli($p$) distribution and red/green signals appear in
contiguous blocks of length $\ell$. \ Let%
\[%
\begin{array}
[c]{c}%
M_{n}=\max\limits_{0\leq j\leq n}S_{j}%
\end{array}
\]
denote the maximum of the walk over the time interval $[0,n]$. \ The quantity
$M_{n}$ is indisputably the most prominent feature of a walk, yet its
first/second moments seem not to be widely known.

A method for generating exact $M_{n}$ probabilities, developed in
\cite{PP-trffc} for symmetric simple RRWs, was extended in \cite{F1-trffc} to
asymmetric simple RRWs. \ This is not to be confused with procedures for
generating exact $S_{n}$ probabilities, which already have a substantial
literature \cite{T1-trffc, T2-trffc, T3-trffc, T4-trffc, T5-trffc, T6-trffc,
T7-trffc, T8-trffc, T9-trffc}.

For the special case $\ell=1$, we have a recurrence for $\mathbb{P}\left\{
S_{n}=x\text{ and }M_{n}=a\right\}  $:%

\[
F_{1}(x,a)=p\,\delta_{x,1}\delta_{a,1}+q\,\delta_{x,0}\delta_{a,0},
\]%
\[
F_{n}(x,a)=\left\{
\begin{array}
[c]{lll}%
p\,\delta_{x,a}\,F_{n-1}(x-1,a-1)+p\,F_{n-1}(x-1,a)+q\,F_{n-1}(x,a) &  &
\text{if }n\text{ is odd,}\\
(p+q\,\delta_{x,0})\,F_{n-1}(x,a)+q(1-\delta_{x,a})F_{n-1}(x+1,a) &  &
\text{if }n\text{ is even}%
\end{array}
\right.
\]
where $x=0,1,\ldots,a$. \ Summing over all $x$ gives $\mathbb{P}\left\{
M_{n}=a\right\}  $. \ It is natural to parse the sequence $\{S_{j}\}$ into
blocks of length two, corresponding to an $RG$ light cycle. \ Another way of
generating exact $M_{2n}$ probabilities\ is to determine the coefficient of
$\lambda^{n}$ in the series expansion of%
\[
\frac{\theta-2\,p\,q}{1-\lambda}\,\left[  \frac{2^{a}p^{2a-1}\theta^{a}%
}{2^{2a}p^{2a-1}q^{2a+1}\left(  \theta-2p^{2}\right)  +\theta^{2a}\left(
\theta-2q^{2}\right)  }-\frac{2^{a+1}p^{2a+1}\theta^{a+1}}{2^{2a+2}%
p^{2a+1}q^{2a+3}\left(  \theta-2p^{2}\right)  +\theta^{2a+2}\left(
\theta-2q^{2}\right)  }\right]
\]
where%
\[
\theta(\lambda)=\dfrac{1-2\,p\,q\,\lambda-\sqrt{1-4\,p\,q\,\lambda}}{\lambda
}.
\]

For $\ell=2$, we have an identical recurrence except the two conditions
\textquotedblleft$n$ is odd\textquotedblright\ versus \textquotedblleft$n$ is
even\textquotedblright\ are replaced by \textquotedblleft$n\equiv
1,2\operatorname{mod}4$\textquotedblright\ versus \textquotedblleft%
$n\equiv3,4\operatorname{mod}4$\textquotedblright\ respectively. \ It is
natural here to parse the sequence $\{S_{j}\}$ into blocks of length four,
corresponding to an $RRGG$ light cycle. \ An expression for $M_{4n}$ analogous
to before, however, seems not to be possible.

\section{Algebra for $\ell=1$}

Combining elements of the recurrence for $\ell=1$ when $n$ is even, we have%
\begin{align*}
F_{n}(x,a)  & =(p+q\,\delta_{x,0})\,F_{n-1}(x,a)+q(1-\delta_{x,a}%
)F_{n-1}(x+1,a)\\
& =(p+q\,\delta_{x,0})\left[  p\,\delta_{x,a}\,F_{n-2}(x-1,a-1)+p\,F_{n-2}%
(x-1,a)+q\,F_{n-2}(x,a)\right]  +\\
& q(1-\delta_{x,a})\left[  p\,\delta_{x+1,a}\,F_{n-2}(x,a-1)+p\,F_{n-2}%
(x,a)+q\,F_{n-2}(x+1,a)\right] \\
& =p\,\delta_{x,a}(p+q\,\delta_{x,0})F_{n-2}(x-1,a-1)+p(p+q\,\delta
_{x,0})F_{n-2}(x-1,a)+\\
& p\,q\,\delta_{x+1,a}(1-\delta_{x,a})F_{n-2}(x,a-1)+\\
& q\left[  (p+q\,\delta_{x,0})+p(1-\delta_{x,a})\right]  F_{n-2}%
(x,a)+q^{2}(1-\delta_{x,a})F_{n-2}(x+1,a)
\end{align*}
for all $n\geq1$, $0\leq x\leq a$ and $a\geq1$. \ Our goal is to solve for
$G_{n}(x,a)=F_{2n}(x,a)$. \ Introduce the generating function%
\[
G(\lambda,x,a)=%
{\displaystyle\sum\limits_{n=1}^{\infty}}
\lambda^{n}G_{n}(x,a)=\lambda\,G_{1}(x,a)+%
{\displaystyle\sum\limits_{n=1}^{\infty}}
\lambda^{n+1}G_{n+1}(x,a)
\]
from which%
\begin{align*}
G(\lambda,x,a)  & =\lambda\,G_{1}(x,a)+\lambda\,p\,\delta_{x,a}(p+q\,\delta
_{x,0})G(\lambda,x-1,a-1)+\lambda\,p(p+q\,\delta_{x,0})G(\lambda,x-1,a)+\\
& \lambda\,p\,q\,\delta_{x+1,a}(1-\delta_{x,a})G(\lambda,x,a-1)+\lambda
\,q\left[  (p+q\,\delta_{x,0})+p(1-\delta_{x,a})\right]  G(\lambda,x,a)+\\
& \lambda\,q^{2}(1-\delta_{x,a})G(\lambda,x+1,a)
\end{align*}
follows. Introduce the double generating function\
\[
\tilde{G}(\lambda,\mu,a)=%
{\displaystyle\sum\limits_{n=1}^{\infty}}
{\displaystyle\sum\limits_{x=0}^{a}}
\lambda^{n}\mu^{x}G_{n}(x,a)=%
{\displaystyle\sum\limits_{x=0}^{a}}
G(\lambda,x,a)\mu^{x}
\]
and note that because%
\[
G_{1}(x,a)=p^{2}\delta_{x,1}\delta_{a,1}+p\,q\,\delta_{x,0}\delta
_{a,1}+q\,\delta_{x,0}\delta_{a,0}
\]
we have%
\[
\lambda%
{\displaystyle\sum\limits_{x=0}^{a}}
G_{1}(x,a)\mu^{x}=\lambda\,q\,\delta_{a,0}+\lambda\,p\,q\,\delta_{a,1}%
+\lambda\,p^{2}\mu\,\delta_{a,1}.
\]
Also,%
\[
\lambda\,p^{2}%
{\displaystyle\sum\limits_{x=0}^{a}}
\delta_{x,a}G(\lambda,x-1,a-1)\mu^{x}+\lambda\,p\,q%
{\displaystyle\sum\limits_{x=0}^{a}}
\delta_{x,0}\delta_{x,a}G(\lambda,x-1,a-1)\mu^{x}=\lambda\,p^{2}%
G(\lambda,a-1,a-1)\mu^{a}+0,
\]%
\[
\lambda\,p^{2}%
{\displaystyle\sum\limits_{x=0}^{a}}
G(\lambda,x-1,a)\mu^{x}+\lambda\,p\,q%
{\displaystyle\sum\limits_{x=0}^{a}}
\delta_{x,0}G(\lambda,x-1,a)\mu^{x}=\lambda\,p^{2}\mu\left(  \tilde{G}%
(\lambda,\mu,a)-G(\lambda,a,a)\mu^{a}\right)  +0,
\]%
\[
\lambda\,p\,q%
{\displaystyle\sum\limits_{x=0}^{a}}
\delta_{x+1,a}(1-\delta_{x,a})G(\lambda,x,a-1)\mu^{x}=\lambda\,p\,q\,G(\lambda
,a-1,a-1)\mu^{a-1},
\]%
\[
\lambda\,p\,q%
{\displaystyle\sum\limits_{x=0}^{a}}
G(\lambda,x,a)\mu^{x}+\lambda\,q^{2}%
{\displaystyle\sum\limits_{x=0}^{a}}
\delta_{x,0}G(\lambda,x,a)\mu^{x}=\lambda\,p\,q\,\tilde{G}(\lambda
,\mu,a)+\lambda\,q^{2}G(\lambda,0,a)\mu^{0},
\]%
\[
\lambda\,p\,q%
{\displaystyle\sum\limits_{x=0}^{a}}
(1-\delta_{x,a})G(\lambda,x,a)\mu^{x}=\lambda\,p\,q\left(  \tilde{G}%
(\lambda,\mu,a)-G(\lambda,a,a)\mu^{a}\right)  ,
\]%
\[
\lambda\,q^{2}%
{\displaystyle\sum\limits_{x=0}^{a}}
(1-\delta_{x,a})G(\lambda,x+1,a)\mu^{x}=\frac{\lambda\,q^{2}}{\mu}\left(
\tilde{G}(\lambda,\mu,a)-G(\lambda,0,a)\mu^{0}\right)  .
\]
We obtain%
\begin{align*}
\tilde{G}  & =\left(  \lambda\,q\,\delta_{a,0}+\lambda\,p\,q\,\delta
_{a,1}+\lambda\,p^{2}\mu\,\delta_{a,1}\right)  +\lambda\,p^{2}G(\lambda
,a-1,a-1)\mu^{a}\\
& +\lambda\,p^{2}\mu\left(  \tilde{G}(\lambda,\mu,a)-G(\lambda,a,a)\mu
^{a}\right)  +\lambda\,p\,q\,G(\lambda,a-1,a-1)\mu^{a-1}\\
& +\lambda\,p\,q\,\tilde{G}(\lambda,\mu,a)+\lambda\,q^{2}G(\lambda
,0,a)+\lambda\,p\,q\left(  \tilde{G}(\lambda,\mu,a)-G(\lambda,a,a)\mu
^{a}\right) \\
& +\frac{\lambda\,q^{2}}{\mu}\left(  \tilde{G}(\lambda,\mu,a)-G(\lambda
,0,a)\right)
\end{align*}
that is,%
\begin{align*}
\left(  1-\lambda\,p^{2}\mu-2\lambda\,p\,q-\frac{\lambda\,q^{2}}{\mu}\right)
\tilde{G}  & =\left(  \lambda\,q\,\delta_{a,0}+\lambda\,p\,q\,\delta
_{a,1}+\lambda\,p^{2}\mu\,\delta_{a,1}\right) \\
& +\left(  \lambda\,q^{2}-\frac{\lambda\,q^{2}}{\mu}\right)  G(\lambda,0,a)\\
& -\left(  \lambda\,p^{2}\mu^{a+1}+\lambda\,p\,q\,\mu^{a}\right)
G(\lambda,a,a)\\
& +\left(  \lambda\,p^{2}\mu^{a}+\lambda\,p\,q\,\mu^{a-1}\right)
G(\lambda,a-1,a-1)
\end{align*}
that is,%
\begin{align}
\left(  p^{2}\mu^{2}+2\,p\,q\,\mu-\frac{\mu}{\lambda}+q^{2}\right)  \tilde{G}
& =-q\,\mu\,\delta_{a,0}-p\,q\,\mu\,\delta_{a,1}-p^{2}\mu^{2}\delta_{a,1}\\
& -q^{2}\left(  \mu-1\right)  G(\lambda,0,a)\nonumber\\
& +p\left(  p\,\mu^{a+2}+q\,\mu^{a+1}\right)  G(\lambda,a,a)\nonumber\\
& -p\left(  \,p\,\mu^{a+1}+q\,\mu^{a}\right)  G(\lambda,a-1,a-1)\nonumber
\end{align}
after multiplying both sides by $-\mu/\lambda$. \ Examine the special case
$a=1$:%
\begin{align*}
p^{2}\left(  \mu-\frac{\theta}{2p^{2}}\right)  \left(  \mu-\frac{2q^{2}%
}{\theta}\right)  \tilde{G}  & =-p\,q\,\mu-p^{2}\mu^{2}-q^{2}\left(
\mu-1\right)  G(\lambda,0,1)\\
& +p\left(  p\,\mu^{3}+q\,\mu^{2}\right)  G(\lambda,1,1)-p\left(  \,p\,\mu
^{2}+q\,\mu\right)  G(\lambda,0,0)
\end{align*}
where%
\[%
\begin{array}
[c]{ccccc}%
\dfrac{\theta}{2p^{2}}=\dfrac{1-2\,p\,q\,\lambda-\sqrt{1-4\,p\,q\,\lambda}%
}{2p^{2}\lambda}, &  & \text{equivalently,} &  & \dfrac{2q^{2}}{\theta}%
=\dfrac{1-2\,p\,q\,\lambda+\sqrt{1-4\,p\,q\,\lambda}}{2p^{2}\lambda}.
\end{array}
\]
For future reference, the sum of the two zeroes is $1/(p^{2}\lambda)-2q/p$,
which implies that
\[
\lambda=\frac{1}{p^{2}}\frac{1}{\frac{2q}{p}+\tfrac{\theta}{2p^{2}}%
+\tfrac{2q^{2}}{\theta}}=\frac{2\,\theta}{\left(  \theta+2\,p\,q\right)  ^{2}%
},
\]%
\[
\frac{p\,\lambda}{1-q\,\lambda}=\frac{\frac{2\,p\,\theta}{\left(
\theta+2\,p\,q\right)  ^{2}}}{1-\frac{2\,q\,\theta}{\left(  \theta
+2\,p\,q\right)  ^{2}}}=\frac{2\,p\,\theta}{\theta^{2}-2(1-2\,p)q\,\theta
+4p^{2}q^{2}}.
\]
Taking $\mu=\theta/(2p^{2})$ and then $\mu=2q^{2}/\theta$, we have%
\[
\left\{
\begin{array}
[c]{l}%
0=-\dfrac{q\,\theta}{2p}-\dfrac{\theta^{2}}{4p^{2}}-q^{2}\left(  \dfrac
{\theta}{2p^{2}}-1\right)  G(\lambda,0,1)+p\left(  \dfrac{\theta^{3}}{8p^{5}%
}+\dfrac{q\,\theta^{2}}{4p^{4}}\right)  G(\lambda,1,1)\\
\;\;\;\;\;\;\,-p\left(  \dfrac{\theta^{2}}{4p^{3}}+\dfrac{q\,\theta}{2p^{2}%
}\right)  G(\lambda,0,0)\\
0=-\dfrac{2\,p\,q^{3}}{\theta}-\dfrac{4p^{2}q^{4}}{\theta^{2}}-q^{2}\left(
\dfrac{2q^{2}}{\theta}-1\right)  G(\lambda,0,1)+p\left(  \dfrac{8\,p\,q^{6}%
}{\theta^{3}}+\dfrac{4q^{5}}{\theta^{2}}\right)  G(\lambda,1,1)\\
\;\;\;\;\;\;\,-p\left(  \dfrac{4\,p\,q^{4}}{\theta^{2}}+\dfrac{2q^{3}}{\theta
}\right)  G(\lambda,0,0)
\end{array}
\right.
\]
and%
\[
G(\lambda,0,0)=%
{\displaystyle\sum\limits_{n=1}^{\infty}}
\lambda^{n}G_{n}(0,0)=%
{\displaystyle\sum\limits_{n=1}^{\infty}}
\lambda^{n}q^{n}=\frac{q\,\lambda\,}{1-q\,\lambda}=\frac{2\,q\,\theta
}{4\,p\,q(\theta+p\,q)+\theta(\theta-2q)};
\]
thus, on eliminating $G(\lambda,0,1)$,
\[
G(\lambda,1,1)=\frac{1}{1-q\,\lambda}\,\frac{2p^{2}\left[  \theta
^{2}-2(1-2\,p)q\,\theta+4p^{2}q^{2}\right]  \left(  \theta-2\,p\,q\right)
\theta}{2^{4}p^{3}q^{5}\left(  \theta-2p^{2}\right)  +\theta^{4}\left(
\theta-2q^{2}\right)  }.
\]
Now examine the general case $a>1$:%
\[
\left\{
\begin{array}
[c]{l}%
0=-q^{2}\left(  \dfrac{\theta}{2p^{2}}-1\right)  G(\lambda,0,a)+p\left(
\dfrac{\theta^{a+2}}{2^{a+2}p^{2a+3}}+\dfrac{q\,\theta^{a+1}}{2^{a+1}p^{2a+2}%
}\right)  G(\lambda,a,a)\\
\;\;\;\;\;\;\;-p\left(  \dfrac{\theta^{a+1}}{2^{a+1}p^{2a+1}}+\dfrac
{q\,\theta^{a}}{2^{a}p^{2a}}\right)  G(\lambda,a-1,a-1)\\
0=-q^{2}\left(  \dfrac{2q^{2}}{\theta}-1\right)  G(\lambda,0,a)+p\left(
\dfrac{2^{a+2}p\,q^{2a+4}}{\theta^{a+2}}+\dfrac{2^{a+1}q^{2a+3}}{\theta^{a+1}%
}\right)  G(\lambda,a,a)\\
\;\;\;\;\;\;\;-p\left(  \dfrac{2^{a+1}p\,q^{2a+2}}{\theta^{a+1}}+\dfrac
{2^{a}q^{2a+1}}{\theta^{a}}\right)  G(\lambda,a-1,a-1)
\end{array}
\right.
\]
and, on eliminating $G(\lambda,0,a)$,%
\begin{align*}
G(\lambda,a,a)  & =2p^{2}\frac{2^{2a}p^{2a-1}q^{2a+1}\left(  \theta
-2p^{2}\right)  +\theta^{2a}\left(  \theta-2q^{2}\right)  }{2^{2a+2}%
p^{2a+1}q^{2a+3}\left(  \theta-2p^{2}\right)  +\theta^{2a+2}\left(
\theta-2q^{2}\right)  }\theta\,G(\lambda,a-1,a-1)\\
& =\frac{1}{1-q\,\lambda}\,\frac{2^{a}p^{2a}\left[  \theta^{2}%
-2(1-2\,p)q\,\theta+4p^{2}q^{2}\right]  \left(  \theta-2\,p\,q\right)
\theta^{a}}{2^{2a+2}p^{2a+1}q^{2a+3}\left(  \theta-2p^{2}\right)
+\theta^{2a+2}\left(  \theta-2q^{2}\right)  }%
\end{align*}
after iteration. \ Finally, given $a>1$ and taking the limit in formula (1) as
$\mu\rightarrow1$, we have%
\begin{align*}
-\frac{1-\lambda}{\lambda}\tilde{G}  & =\left(  1-\frac{1}{\lambda}\right)
\tilde{G}=\left(  (p+q)^{2}-\frac{1}{\lambda}\right)  \tilde{G}=\left(
p^{2}+2\,p\,q-\frac{1}{\lambda}+q^{2}\right)  \tilde{G}=\\
& =p(p+q)G(\lambda,a,a)-p(p+q)G(\lambda,a-1,a-1)=p\,G(\lambda
,a,a)-p\,G(\lambda,a-1,a-1)
\end{align*}
therefore%
\[
\tilde{G}=\frac{p\,\lambda}{1-\lambda}\left[  G(\lambda,a-1,a-1)-G(\lambda
,a,a)\right]  ,
\]
as was to be shown. \ The case $a=1$ must be treated separately:%
\begin{align*}
\tilde{G}  & =G(\lambda,0,1)+G(\lambda,1,1)\\
& =\frac{2p^{2}}{2p^{2}-\theta}\frac{1}{q^{2}}\left[  \dfrac{q\,\theta}%
{2p}+\dfrac{\theta^{2}}{4p^{2}}-p\left(  \dfrac{\theta^{3}}{8p^{5}}%
+\dfrac{q\,\theta^{2}}{4p^{4}}\right)  G(\lambda,1,1)+p\left(  \dfrac
{\theta^{2}}{4p^{3}}+\dfrac{q\,\theta}{2p^{2}}\right)  G(\lambda,0,0)\right]
+G(\lambda,1,1)\\
& =\frac{p\,\lambda\left(  1-p\,q\,\lambda\right)  }{(1-q\,\lambda)\left[
p\,q^{2}\lambda^{2}-(1+2p)q\,\lambda+1\right]  }%
\end{align*}
consistent with the series expansion in the introduction.

\section{Calculus for $\ell=1$}

Setting $p=q=1/2$, we obtain%
\begin{align*}
\tilde{G}(\lambda,1,a)  & =\frac{1}{1-\lambda}\left[  \frac{\left(  \frac
{1}{2}\right)  ^{a-1}\theta^{a}}{\left(  \frac{1}{2}\right)  ^{2a}+\theta
^{2a}}-\frac{\left(  \frac{1}{2}\right)  ^{a}\theta^{a+1}}{\left(  \frac{1}%
{2}\right)  ^{2a+2}+\theta^{2a+2}}\right]  \\
& =\frac{2}{1-\lambda}\left[  \frac{(2\theta)^{a}}{1+(2\theta)^{2a}}%
-\frac{(2\theta)^{a+1}}{1+(2\theta)^{2a+2}}\right]  \\
& =\frac{2}{1-\lambda}\left[  \frac{1}{(2\theta)^{a}+(2\theta)^{-a}}-\frac
{1}{(2\theta)^{a+1}+(2\theta)^{-a-1}}\right]
\end{align*}
and thus have the double generating function%
\begin{align*}
\Psi(\lambda,\nu)  & =%
{\displaystyle\sum\limits_{n=1}^{\infty}}
{\displaystyle\sum\limits_{a=1}^{\infty}}
\lambda^{n}\nu^{a}\mathbb{P}\left\{  M_{2n}=a\right\}  =%
{\displaystyle\sum\limits_{a=1}^{\infty}}
\nu^{a}\tilde{G}(\lambda,1,a)\\
& =\frac{1}{1-\lambda}%
{\displaystyle\sum\limits_{a=1}^{\infty}}
\nu^{a}\left[  \frac{2}{(2\theta)^{a}+(2\theta)^{-a}}-\frac{2}{(2\theta
)^{a+1}+(2\theta)^{-a-1}}\right]  \\
& =\frac{1}{1-\lambda}\left[  \nu^{0}\frac{2}{(2\theta)^{1}+(2\theta)^{-1}}+%
{\displaystyle\sum\limits_{a=1}^{\infty}}
\left(  \nu^{a}-\nu^{a-1}\right)  \frac{2}{(2\theta)^{a}+(2\theta)^{-a}%
}\right]  \\
& =\frac{\lambda}{(1-\lambda)(2-\lambda)}-\frac{1-\nu}{1-\lambda}%
{\displaystyle\sum\limits_{a=1}^{\infty}}
\nu^{a-1}\frac{2}{(2\theta)^{a}+(2\theta)^{-a}}.
\end{align*}
Let us focus on $\mathbb{E}\left(  M_{2n}\right)  $:%
\begin{align*}%
{\displaystyle\sum\limits_{n=1}^{\infty}}
\lambda^{n}\mathbb{E}\left(  M_{2n}\right)    & =\left.  \frac{\partial\Psi
}{\partial\nu}\right\vert _{\nu=1}\\
& =\left.  \frac{1}{1-\lambda}%
{\displaystyle\sum\limits_{a=1}^{\infty}}
\nu^{a-1}\frac{2}{(2\theta)^{a}+(2\theta)^{-a}}-\frac{1-\nu}{1-\lambda}%
{\displaystyle\sum\limits_{a=1}^{\infty}}
(a-1)\nu^{a-2}\frac{2}{(2\theta)^{a}+(2\theta)^{-a}}\right\vert _{\nu=1}\\
& =\frac{1}{1-\lambda}%
{\displaystyle\sum\limits_{a=1}^{\infty}}
\frac{2}{(2\theta)^{a}+(2\theta)^{-a}}%
\end{align*}
which provides that (in an extended sense) the following sequence is Abel
convergent \cite{PP-trffc, PS-trffc}:%
\begin{align*}
\lim_{n\rightarrow\infty}^{\;\;\;\;\;\;\;\;\ast}\frac{\mathbb{E}\left(
M_{2n}\right)  }{\sqrt{n}}  & =\lim_{\lambda\rightarrow1^{-}}(1-\lambda)^{3/2}%
{\displaystyle\sum\limits_{n=1}^{\infty}}
\frac{1}{(1/2)!}\lambda^{n-1/2}\mathbb{E}\left(  M_{2n}\right)  \\
& =\lim_{\lambda\rightarrow1^{-}}\left(  \frac{1-\lambda}{\lambda}\right)
^{1/2}%
{\displaystyle\sum\limits_{a=1}^{\infty}}
\frac{2}{(2\theta)^{a}+(2\theta)^{-a}}\frac{1}{(1/2)!}.
\end{align*}
Let $2\theta=\exp(-t)$, then%
\[
\frac{2}{(2\theta)^{a}+(2\theta)^{-a}}=\frac{2}{e^{at}+e^{-at}}%
=\operatorname{sech}(at)
\]
and, because $\lambda=2\theta/\left(  \theta+\frac{1}{2}\right)  ^{2}%
=2\theta/\left(  \theta^{2}+\theta+\frac{1}{4}\right)  $,
\[
\dfrac{1-\lambda}{\lambda}=\frac{\theta^{2}+\frac{1}{4}+\theta}{2\theta
}-1=\frac{(2\theta)+(2\theta)^{-1}}{4}-\frac{1}{2}=\frac{\cosh(t)-1}{2}.
\]
We have%
\begin{align*}
\lim_{n\rightarrow\infty}^{\;\;\;\;\;\;\;\;\ast}\frac{\mathbb{E}\left(
M_{2n}\right)  }{\sqrt{n}}  & =\lim_{t\rightarrow0^{+}}\sqrt{\frac{\cosh
(t)-1}{2}}\frac{1}{(1/2)!}%
{\displaystyle\sum\limits_{a=1}^{\infty}}
\operatorname{sech}(at)\\
& =\lim_{t\rightarrow0^{+}}\sqrt{\frac{1}{\pi}}t%
{\displaystyle\sum\limits_{a=1}^{\infty}}
\operatorname{sech}(at)
\end{align*}
since $(\cosh(t)-1)/2\sim t^{2}/4$ and $(1/2)!=\sqrt{\pi}/2$. \ By a Riemann
sum-based argument,
\[
t%
{\displaystyle\sum\limits_{a=\left\lceil \alpha/t\right\rceil }^{\left\lfloor
\beta/t\right\rfloor }}
\operatorname{sech}(at)\rightarrow%
{\displaystyle\int\limits_{\alpha}^{\beta}}
\operatorname{sech}(b)db
\]
as $t\rightarrow0^{+}$, which in turn gives%
\[
\lim_{n\rightarrow\infty}^{\;\;\;\;\;\;\;\;\ast}\frac{\mathbb{E}\left(
M_{2n}\right)  }{\sqrt{2n}}=\sqrt{\frac{1}{2\pi}}%
{\displaystyle\int\limits_{0}^{\infty}}
\operatorname{sech}(b)db=\sqrt{\frac{\pi}{8}}%
\]
as $\alpha\rightarrow0^{+}$ and $\beta\rightarrow\infty$. \ A similar argument
gives%
\[
\lim_{n\rightarrow\infty}^{\;\;\;\;\;\;\;\;\ast}\frac{\mathbb{E}\left(
M_{2n}^{2}\right)  }{2n}=\frac{1}{4}%
{\displaystyle\int\limits_{0}^{\infty}}
b\operatorname{sech}(b)db=\frac{G}{2}%
\]
where $G$ is Catalan's constant \cite{F2-trffc}.

\section{Algebra for $\ell=2$}

Combining elements of the recurrence for $\ell=2$ when $n\equiv
4\operatorname{mod}4 $, we have%
\begin{align*}
F_{n}(x,a)  & =(p+q\,\delta_{x,0})\,F_{n-1}(x,a)+q(1-\delta_{x,a}%
)F_{n-1}(x+1,a)\\
& =(p+q\,\delta_{x,0})\left[  (p+q\,\delta_{x,0})\,F_{n-2}(x,a)+q(1-\delta
_{x,a})F_{n-2}(x+1,a)\right]  +\\
& q(1-\delta_{x,a})\left[  (p+q\,\delta_{x+1,0})\,F_{n-2}(x+1,a)+q(1-\delta
_{x+1,a})F_{n-2}(x+2,a)\right] \\
& =(p+q\,\delta_{x,0})^{2}F_{n-2}(x,a)+q(2p+q\,\delta_{x,0})(1-\delta
_{x,a})F_{n-2}(x+1,a)+\\
& q^{2}(1-\delta_{x,a})(1-\delta_{x+1,a})F_{n-2}(x+2,a)
\end{align*}
owing to $\delta_{x+1,0}=0$ for all $x\geq0$, and%
\begin{align*}
F_{n-2}(x,a)  & =p\,\delta_{x,a}\,F_{n-3}(x-1,a-1)+p\,F_{n-3}%
(x-1,a)+q\,F_{n-3}(x,a)\\
& =p\,\delta_{x,a}\left[  p\,\delta_{x,a}\,F_{n-4}(x-2,a-2)+p\,F_{n-4}%
(x-2,a-1)+q\,F_{n-4}(x-1,a-1)\right]  +\\
& p\left[  p\,\delta_{x-1,a}\,F_{n-4}(x-2,a-1)+p\,F_{n-4}(x-2,a)+q\,F_{n-4}%
(x-1,a)\right]  +\\
& q\left[  p\,\delta_{x,a}\,F_{n-4}(x-1,a-1)+p\,F_{n-4}(x-1,a)+q\,F_{n-4}%
(x,a)\right] \\
& =p^{2}\delta_{x,a}F_{n-4}(x-2,a-2)+p^{2}\left(  \delta_{x,a}+\delta
_{x-1,a}\right)  F_{n-4}(x-2,a-1)+p^{2}F_{n-4}(x-2,a)+\\
& 2\,p\,q\,\delta_{x,a}F_{n-4}(x-1,a-1)+2\,p\,q\,F_{n-4}(x-1,a)+q^{2}%
F_{n-4}(x,a)
\end{align*}
owing to $\delta_{x-1,a-1}=\delta_{x,a}$ always. \ Let\ $G_{n}(x,a)=F_{4n}%
(x,a)$. \ A recursion for $G_{n}(x,a)$ in terms of%
\begin{align*}
& G_{n-1}(x-2,a-2),\;G_{n-1}(x-1,a-2),\;G_{n-1}(x,a-2),\\
& G_{n-1}(x-2,a-1),\;G_{n-1}(x-1,a-1),\;G_{n-1}(x,a-1),\;G_{n-1}(x+1,a-1),\\
& G_{n-1}(x-2,a),\;G_{n-1}(x-1,a),\;G_{n-1}(x,a),\;G_{n-1}(x+1,a),\;G_{n-1}%
(x+2,a)
\end{align*}
arises (too complicated to reproduce here). \ Note, for example,%
\[
G_{1}(x,a)=p^{4}\delta_{x,2}\delta_{a,2}+2\,p^{3}q\,\delta_{x,1}\delta
_{a,2}+p^{2}q^{2}\delta_{x,0}\delta_{a,2}+2\,p^{3}q\,\delta_{x,1}\delta
_{a,1}+2\,p(1+p)q^{2}\delta_{x,0}\delta_{a,1}+q^{2}\delta_{x,0}\delta_{a,0}
\]
and hence%
\begin{align*}
\lambda%
{\displaystyle\sum\limits_{x=0}^{a}}
G_{1}(x,a)\mu^{x}  & =\lambda\,q^{2}\delta_{a,0}+2\,\lambda\,p(1+p)q^{2}%
\delta_{a,1}+2\,\lambda\,p^{3}q\,\mu\,\delta_{a,1}+\\
& \lambda\,p^{2}q^{2}\delta_{a,2}+2\,\lambda\,p^{3}q\,\mu\,\delta
_{a,2}+\lambda\,p^{4}\mu^{2}\,\delta_{a,2}.
\end{align*}
Recall the quadratic coefficient for $\tilde{G}$ we found in formula (1);
here, this becomes a quartic:%
\begin{equation}
p^{4}\mu^{4}+4\,p^{3}q\,\mu^{3}+6\,p^{2}q^{2}\mu^{2}-\frac{\mu^{2}}{\lambda
}+4\,p\,q^{3}\mu+q^{4}%
\end{equation}
which factors as%
\[
p^{4}\left(  \mu-\frac{\theta}{2p^{2}}\right)  \left(  \mu-\frac{2q^{2}%
}{\theta}\right)  \left(  \mu-\frac{\omega}{2p^{2}}\right)  \left(  \mu
-\frac{2q^{2}}{\omega}\right)
\]
where
\[%
\begin{array}
[c]{ccccc}%
\dfrac{\theta}{2p^{2}}=\dfrac{1-2\,p\,q\sqrt{\lambda}-\sqrt{1-4\,p\,q\sqrt
{\lambda}}}{2p^{2}\sqrt{\lambda}}, &  & \text{equivalently,} &  &
\dfrac{2q^{2}}{\theta}=\dfrac{1-2\,p\,q\sqrt{\lambda}+\sqrt{1-4\,p\,q\sqrt
{\lambda}}}{2p^{2}\sqrt{\lambda}};
\end{array}
\]%
\[%
\begin{array}
[c]{ccccc}%
\dfrac{\omega}{2p^{2}}=\dfrac{-1-2\,p\,q\sqrt{\lambda}-\sqrt{1+4\,p\,q\sqrt
{\lambda}}}{2p^{2}\sqrt{\lambda}}, &  & \text{equivalently,} &  &
\dfrac{2q^{2}}{\omega}=\dfrac{-1-2\,p\,q\sqrt{\lambda}+\sqrt{1+4\,p\,q\sqrt
{\lambda}}}{2p^{2}\sqrt{\lambda}}.
\end{array}
\]
Define%
\begin{align*}
\Delta_{1}  & =256\,p^{8}q^{8}+1024\,p^{7}q^{7}\theta+64\,p^{4}q^{6}%
(22\,p^{2}-q-3\,p\,q)\theta^{2}+\\
& 128\,p^{3}q^{5}(8\,p^{2}-q-3\,p\,q)\theta^{3}+16\,p^{2}q^{4}(35\,p^{2}%
-5\,q-17\,p\,q)\theta^{4}+\\
& 32\,p\,q^{3}(8\,p^{2}-q-3\,p\,q)\theta^{5}+4\,q^{2}(22\,p^{2}%
-q-3\,p\,q)\theta^{6}+16\,p\,q\,\theta^{7}+\theta^{8},
\end{align*}
then it can be shown that%
\[
G(\lambda,0,1)=\frac{8\,p(1+p)q^{2}\left(  \theta+2\,p\,q\right)  ^{4}%
\theta^{2}}{(1-q^{2}\lambda)\Delta_{1}},
\]%
\[
G(\lambda,1,1)=\frac{8\,p^{3}q\left[  \theta^{2}+2(2\,p-1)q\,\theta
+4\,p^{2}q^{2}\right]  \left[  \theta^{2}+2(2\,p+1)q\,\theta+4\,p^{2}%
q^{2}\right]  \theta^{2}}{(1-q^{2}\lambda)\Delta_{1}}.
\]
Unfortunately the expressions become cumbersome beyond this point; no pattern
is evident. \ Define%
\begin{align*}
\Gamma & =\left[  \theta^{2}+2(2p-1)q\,\theta+4p^{2}q^{2}\right]  \left[
\theta^{2}+2(2p+1)q\,\theta+4p^{2}q^{2}\right]  \cdot\\
& \left[  \theta^{4}+8pq\theta^{3}+36p^{2}q^{2}\theta^{2}+32p^{3}q^{3}%
\theta+16p^{4}q^{4}\right]  \cdot\\
& \left\{  4p^{2}q^{2}(\theta+2pq)^{4}+\left[  4pq(\theta+pq)(\theta
+4pq)\left(  \theta^{2}-2q^{2}\theta+4p^{2}q^{2}\right)  \right]
\omega+\right. \\
& \left[  \theta^{4}+2(3p-1)q\,\theta^{3}-8p(1+q)q^{2}\theta^{2}%
+8p^{2}(3p-1)q^{3}\theta+16p^{4}q^{4}\right]  \omega^{2}-\\
& \left.  2(1+p)q\,\theta^{2}\omega^{3}\right\}  ,
\end{align*}%
\begin{align*}
\Delta_{2}  & =4096p^{12}q^{12}+8192p^{11}q^{11}\theta-2048p^{10}q^{12}%
\theta+7168p^{10}q^{10}\theta^{2}-4096p^{9}q^{11}\theta^{2}+5120p^{9}%
q^{9}\theta^{3}-\\
& 3072p^{8}q^{10}\theta^{3}+3072p^{8}q^{8}\theta^{4}-1536p^{7}q^{9}\theta
^{4}+1280p^{7}q^{7}\theta^{5}-768p^{6}q^{8}\theta^{5}+448p^{6}q^{6}\theta
^{6}-\\
& 256p^{5}q^{7}\theta^{6}+128p^{5}q^{5}\theta^{7}-32p^{4}q^{6}\theta
^{7}+16p^{4}q^{4}\theta^{8}-2048p^{10}q^{12}\omega+2048p^{10}q^{10}%
\theta\omega+\\
& 1024p^{8}q^{11}\theta\omega-5120p^{9}q^{11}\theta\omega+6144p^{9}q^{9}%
\theta^{2}\omega+2048p^{7}q^{10}\theta^{2}\omega-5632p^{8}q^{10}\theta
^{2}\omega+\\
& 7424p^{8}q^{8}\theta^{3}\omega+1536p^{6}q^{9}\theta^{3}\omega-4096p^{7}%
q^{9}\theta^{3}\omega+4864p^{7}q^{7}\theta^{4}\omega+768p^{5}q^{8}\theta
^{4}\omega-\\
& 2304p^{6}q^{8}\theta^{4}\omega+1856p^{6}q^{6}\theta^{5}\omega+384p^{4}%
q^{7}\theta^{5}\omega-1024p^{5}q^{7}\theta^{5}\omega+384p^{5}q^{5}\theta
^{6}\omega+\\
& 128p^{3}q^{6}\theta^{6}\omega-352p^{4}q^{6}\theta^{6}\omega+32p^{4}%
q^{4}\theta^{7}\omega+16p^{2}q^{5}\theta^{7}\omega-80p^{3}q^{5}\theta
^{7}\omega-8p^{2}q^{4}\theta^{8}\omega+\\
& 1024p^{10}q^{10}\omega^{2}+4096p^{9}q^{9}\theta\omega^{2}-512p^{8}%
q^{10}\theta\omega^{2}+6144p^{8}q^{8}\theta^{2}\omega^{2}-2048p^{7}q^{9}%
\theta^{2}\omega^{2}+\\
& 5120p^{7}q^{7}\theta^{3}\omega^{2}-2944p^{6}q^{8}\theta^{3}\omega
^{2}+2944p^{6}q^{6}\theta^{4}\omega^{2}-2048p^{5}q^{7}\theta^{4}\omega
^{2}+1280p^{5}q^{5}\theta^{5}\omega^{2}-\\
& 736p^{4}q^{6}\theta^{5}\omega^{2}+384p^{4}q^{4}\theta^{6}\omega^{2}%
-128p^{3}q^{5}\theta^{6}\omega^{2}+64p^{3}q^{3}\theta^{7}\omega^{2}%
-8p^{2}q^{4}\theta^{7}\omega^{2}+4p^{2}q^{2}\theta^{8}\omega^{2}+\\
& 1024p^{9}q^{9}\omega^{3}+3328p^{8}q^{8}\theta\omega^{3}+3328p^{7}q^{7}%
\theta^{2}\omega^{3}-256p^{5}q^{8}\theta^{2}\omega^{3}-896p^{6}q^{8}\theta
^{2}\omega^{3}+\\
& 1984p^{6}q^{6}\theta^{3}\omega^{3}-256p^{4}q^{7}\theta^{3}\omega
^{3}-1280p^{5}q^{7}\theta^{3}\omega^{3}+1088p^{5}q^{5}\theta^{4}\omega
^{3}-64p^{3}q^{6}\theta^{4}\omega^{3}-\\
& 768p^{4}q^{6}\theta^{4}\omega^{3}+496p^{4}q^{4}\theta^{5}\omega^{3}%
-64p^{2}q^{5}\theta^{5}\omega^{3}-320p^{3}q^{5}\theta^{5}\omega^{3}%
+208p^{3}q^{3}\theta^{6}\omega^{3}-\\
& 16pq^{4}\theta^{6}\omega^{3}-56p^{2}q^{4}\theta^{6}\omega^{3}+52p^{2}%
q^{2}\theta^{7}\omega^{3}+4pq\theta^{8}\omega^{3}+256p^{8}q^{8}\omega
^{4}+768p^{7}q^{7}\theta\omega^{4}+\\
& 640p^{6}q^{6}\theta^{2}\omega^{4}-64p^{4}q^{7}\theta^{2}\omega^{4}%
-192p^{5}q^{7}\theta^{2}\omega^{4}+320p^{5}q^{5}\theta^{3}\omega^{4}%
-64p^{3}q^{6}\theta^{3}\omega^{4}-\\
& 224p^{4}q^{6}\theta^{3}\omega^{4}+176p^{4}q^{4}\theta^{4}\omega^{4}%
-16p^{2}q^{5}\theta^{4}\omega^{4}-112p^{3}q^{5}\theta^{4}\omega^{4}%
+80p^{3}q^{3}\theta^{5}\omega^{4}-\\
& 16pq^{4}\theta^{5}\omega^{4}-56p^{2}q^{4}\theta^{5}\omega^{4}+40p^{2}%
q^{2}\theta^{6}\omega^{4}-4q^{3}\theta^{6}\omega^{4}-12pq^{3}\theta^{6}%
\omega^{4}+12pq\theta^{7}\omega^{4}+\theta^{8}\omega^{4},
\end{align*}
then it can be shown that%
\[
G(\lambda,0,2)=\frac{4p^{2}q^{2}(\theta+2\,p\,q)^{2}\theta\,\omega\,\Gamma
}{(1-q^{2}\lambda)\Delta_{1}\Delta_{2}}.
\]
We omit analogous expressions for $G(\lambda,1,2)$ and $G(\lambda,2,2)$.
\ Computer simulations suggest that first/second moments of $M_{n}/\sqrt{n}$
for $p=q=1/2$ and $\ell\geq2$ are numerically equal to those for $\ell=1$,
given suitably large $n$. \ It is reasonable to conjecture that the maximum
queue length distribution $\digamma$ enjoys a kind of universality, depending
only on the values $p\leq q$. \ The function $\digamma$ is apparently
independent of the choice of light cycles $RG$, $RRGG$, $RRRGGG$, \ldots\ but
a rigorous proof may be difficult.

\section{Asymptotic Distribution of $S_{n}$}

No universality applies to the asymptotic distribution of $S_{n}$. \ The light
cycle $RG$ is distinct from $RRGG$ in this regard. \ For $\ell=1$, we have%
\[%
\begin{array}
[c]{ccc}%
\lim\limits_{n\rightarrow\infty}\mathbb{P}\left\{  S_{n}=x\right\}
=\dfrac{(q-p)p^{2x}}{q^{2x+2}}, &  & x=0,1,2,\ldots
\end{array}
\]
which is found via techniques in \cite{T8-trffc}. \ For $\ell=2$, the
corresponding probability for $x=0$ is%
\[
\frac{4(q-p)}{q^{2}\left(  1+2\,q+\sqrt{1+4\,p\,q}\right)  }
\]
and for $x=1$ is%
\[
\frac{4(q-p)\left[  1+2\,p\,q(q-p)-(q-p+2\,p\,q)\sqrt{1+4\,p\,q}\right]
}{q^{4}\left[  -(q-p)+\sqrt{1+4\,p\,q}\right]  \left(  1+2\,q+\sqrt
{1+4\,p\,q}\right)  }.
\]
Let us elaborate on how these formulas are derived. \ From the equation%
\[
z^{\ell}=(q+p\,z)^{2\,\ell}
\]
we obtain $\ell$ roots%
\[
z_{0}=1,\;z_{1},\;z_{2},\;\ldots,\;z_{\ell-1}
\]
satisfying $\left\vert z_{k}\right\vert \leq1$ for all $k$, and determine
$w_{0},\,w_{1},\,w_{2},\,\ldots,\,w_{\ell-1}$ via the linear system%
\[
\left(
\begin{array}
[c]{ccccc}%
1 & 1 & 1 & \ldots & 1\\
1 & \dfrac{z_{1}}{q+p\,z_{1}} & \left(  \dfrac{z_{1}}{q+p\,z_{1}}\right)  ^{2}
& \ldots & \left(  \dfrac{z_{1}}{q+p\,z_{1}}\right)  ^{\ell-1}\\
1 & \dfrac{z_{2}}{q+p\,z_{2}} & \left(  \dfrac{z_{2}}{q+p\,z_{2}}\right)  ^{2}
& \ldots & \left(  \dfrac{z_{2}}{q+p\,z_{2}}\right)  ^{\ell-1}\\
\vdots & \vdots & \vdots &  & \vdots\\
1 & \dfrac{z_{\ell-1}}{q+p\,z_{\ell-1}} & \left(  \dfrac{z_{\ell-1}%
}{q+p\,z_{\ell-1}}\right)  ^{2} & \ldots & \left(  \dfrac{z_{\ell-1}%
}{q+p\,z_{\ell-1}}\right)  ^{\ell-1}%
\end{array}
\right)  \left(
\begin{array}
[c]{c}%
w_{0}\\
w_{1}\\
w_{2}\\
\vdots\\
w_{\ell-1}%
\end{array}
\right)  =\left(
\begin{array}
[c]{c}%
\dfrac{q-p}{q}\ell\\
0\\
0\\
\vdots\\
0
\end{array}
\right)  .
\]
Then%
\[
H(z)=\frac{(q+p\,z)^{\,\ell}}{z^{\ell}-(q+p\,z)^{2\,\ell}}\left(  \frac
{z}{q+p\,z}-1\right)
{\displaystyle\sum\limits_{k=0}^{\ell-1}}
w_{k}\left(  \dfrac{z}{q+p\,z}\right)  ^{k}
\]
is the probability generating function for $S_{n}$ as $n\rightarrow\infty$.
\ In particular, for $\ell=1$,%
\[%
\begin{array}
[c]{ccc}%
w_{0}=\dfrac{q-p}{q}, &  & H(z)=\dfrac{(q-p)(z-1)}{z-(q+p\,z)^{2}}=\dfrac
{q-p}{q^{2}-p^{2}z},
\end{array}
\]%
\[
\lim_{n\rightarrow\infty}\mathbb{E}\left(  S_{n}\right)  =H^{\prime
}(1)=\left.  \frac{p^{2}\left(  q-p\right)  }{\left(  q^{2}-p^{2}z\right)
^{2}}\right\vert _{z=1}=\frac{p^{2}}{q-p},
\]%
\[
\lim_{n\rightarrow\infty}\mathbb{E}\left(  S_{n}(S_{n}-1)\right)
=H^{^{\prime\prime}}(1)=\left.  \frac{2p^{4}\left(  q-p\right)  }{\left(
q^{2}-p^{2}z\right)  ^{3}}\right\vert _{z=1}=\frac{2p^{4}}{(q-p)^{2}}
\]
and, for $\ell=2$,%
\[%
\begin{array}
[c]{ccc}%
w_{0}=\dfrac{4(q-p)}{2+(q-p)+\sqrt{1+4\,p\,q}}, &  & w_{1}=\dfrac
{4\,p(q-p)}{q\left[  -(q-p)+\sqrt{1+4\,p\,q}\right]  },
\end{array}
\]%
\begin{align*}
H(z)  & =\frac{4\,q(q-p)(q+p\,z)}{\left(  q^{2}-p^{2}z\right)  \left[
q^{2}+(1+2\,p\,q)z+p^{2}z^{2}\right]  }\cdot\\
& \left\{  \frac{1}{2+(q-p)+\sqrt{1+4\,p\,q}}+\dfrac{p\,z}{q(q+p\,z)\left[
-(q-p)+\sqrt{1+4\,p\,q}\right]  }\right\}  ,
\end{align*}%
\[
\lim_{n\rightarrow\infty}\mathbb{E}\left(  S_{n}\right)  =H^{\prime}%
(1)=\frac{1}{4}\left[  -4+2(q-p)+\frac{1}{q-p}+\sqrt{1+4\,p\,q}\right]  .
\]
We obtained $H(z)$ for $\ell=2$ by noting that%
\[
(q+p\,z)^{4}-z^{2}=(1-z)(q^{2}-p^{2}z)\left[  q^{2}+(1+2\,p\,q)z+p^{2}%
z^{2}\right]  ,
\]
simplifying the first factor. \ In words, $\lim_{n\rightarrow\infty}%
\mathbb{E}\left(  S_{n}\right)  $ follows a geometric distribution for $\ell=1
$; \ its mean for $\ell=2$ is noticeably smaller than that for $\ell=1$.

\section{Appendix}

Here is more information pertinent to Section 3. \ For $a=1$, we solve a
$2\times2$ system:
\begin{align*}
0  & =-\left[  2p^{3}q\mu+2p(1+p)q^{2}\right]  \mu^{2}-\left[  2p^{3}%
q\mu+4p^{2}q^{2}+2pq^{3}\right]  \mu^{2}G(\lambda,0,0)+\\
& \left[  p^{4}\mu^{4}+2p^{3}q\mu^{3}+p^{2}q^{2}\mu^{2}-(1+3p)q^{3}\mu
^{2}+4pq^{3}\mu+q^{4}\right]  G(\lambda,0,1)+\\
& \left[  p^{4}\mu^{4}+2p^{3}q\mu^{3}+2p^{3}q\mu^{3}+5p^{2}q^{2}\mu
^{2}+2pq^{3}\mu-q^{4}\mu+q^{4}\right]  \mu\,G(\lambda,1,1)
\end{align*}
for $G(\lambda,0,1)$ and $G(\lambda,1,1)$, taking $\mu=\theta/(2p^{2})$ and
$\mu=2q^{2}/\theta$, utilizing%
\[
G(\lambda,0,0)=%
{\displaystyle\sum\limits_{n=1}^{\infty}}
\lambda^{n}G_{n}(0,0)=%
{\displaystyle\sum\limits_{n=1}^{\infty}}
\lambda^{n}q^{2n}=\frac{q^{2}\lambda\,}{1-q^{2}\lambda}.
\]
For $a=2$, we solve a $3\times3$ system:%
\begin{align*}
0  & =-\left[  p^{4}\mu^{2}+2p^{3}q\mu+p^{2}q^{2}\right]  \mu^{2}-\left[
p^{4}\mu^{2}+2p^{3}q\mu+p^{2}q^{2}\right]  \mu^{2}G(\lambda,0,0)-\\
& \left[  p^{4}\mu^{2}+2p^{3}q\mu+p^{2}q^{2}\right]  \mu^{2}G(\lambda,0,1)-\\
& \left[  2p^{3}q\mu^{2}+4p^{2}q^{2}\mu+2pq^{3}\right]  \mu^{2}G(\lambda
,1,1)+\\
& \left[  p^{4}\mu^{4}+2p^{3}q\mu^{3}+p^{2}q^{2}\mu^{2}-q^{4}\mu+q^{4}\right]
\mu\,G(\lambda,1,2)+\\
& \left[  p^{4}\mu^{3}+4p^{3}q\mu^{2}+5p^{2}q^{2}\mu+2pq^{3}\right]  \mu
^{3}G(\lambda,2,2)-\\
& \left[  (1+3p)q^{3}\mu^{2}-4pq^{3}\mu-q^{4}\right]  G(\lambda,0,2)
\end{align*}
for $G(\lambda,0,2)$, $G(\lambda,1,2)$ and $G(\lambda,2,2)$; taking
$\mu=\theta/(2p^{2})$, $\mu=2q^{2}/\theta$ and $\mu=\omega/(2p^{2})$. \ For
$a\geq3$, we solve a $4\times4$ system:%
\begin{align*}
0  & =-\left[  p^{4}\mu^{2}+2p^{3}q\mu+p^{2}q^{2}\right]  \mu^{a}%
G(\lambda,a-2,a-2)-\\
& \left[  p^{4}\mu^{2}+2p^{3}q\mu+p^{2}q^{2}\right]  \mu^{a}G(\lambda
,a-2,a-1)-\\
& \left[  2p^{3}q\mu^{2}+4p^{2}q^{2}\mu+2pq^{3}\right]  \mu^{a}G(\lambda
,a-1,a-1)+\\
& \left[  p^{4}\mu^{2}+2p^{3}q\mu+p^{2}q^{2}\right]  \mu^{a+1}G(\lambda
,a-1,a)+\\
& \left[  p^{4}\mu^{3}+4p^{3}q\mu^{2}+5p^{2}q^{2}\mu+2pq^{3}\right]  \mu
^{a+1}G(\lambda,a,a)-\\
& \left[  (1+3p)q^{3}\mu^{2}-4pq^{3}\mu-q^{4}\right]  G(\lambda,0,a)-q^{4}%
\left(  \mu-1\right)  \mu\,G(\lambda,1,a)
\end{align*}
for $G(\lambda,0,a)$, $G(\lambda,1,a)$, $G(\lambda,a-1,a)$ and $G(\lambda
,a,a)$; taking $\mu=\theta/(2p^{2})$, $\mu=2q^{2}/\theta$, $\mu=\omega
/(2p^{2})$ and $\mu=2q^{2}/\omega$. \ Note that, via this procedure, the terms
$G(\lambda,2,4)$, $G(\lambda,2,5)$, $G(\lambda,3,5)$, $G(\lambda,2,6)$,
$G(\lambda,3,6)$, $G(\lambda,4,6)$, \ldots\ remain open. \ Formula (2)
experiences only limited success in calculating $\tilde{G}(\lambda,\mu,a)$,
unlike formula (1).

\section{Acknowledgements}

I am thankful to Marko Boon \cite{T9-trffc} for helpful discussions that led
to Section 4.

\end{document}